\documentclass[graybox]{svmult}

\usepackage{mathptmx}       
\usepackage{helvet}         
\usepackage{courier}        
\usepackage{type1cm}        
\usepackage{makeidx}         
\usepackage{graphicx}        
\usepackage{multicol}        
\usepackage[bottom]{footmisc}

\makeindex             

\usepackage{amsmath,amssymb}

\DeclareMathOperator*{\argmin}{{\rm argmin}}

\def\R{\mathbb{R}}
\def\N{\mathbb{N}}
\def\T{\mathbb{T}}
\def\E{\mathbb{E}}
\def\P{\mathbb{P}}
\def\div{{\rm div}}

\begin{document}

\title*{Two mathematical tools to analyze metastable stochastic processes}
\author{Tony Leli\`evre}
\institute{Tony Leli\`evre \at CERMICS, Ecole des Ponts Paristech, Université Paris-Est, 6 et 8 avenue Blaise Pascal, 77455 Marne-la-Vall\'ee, France and INRIA Rocquencourt, MicMac project team, Domaine de Voluceau, B.P. 105, 78153 Le Chesnay Cedex, France. \email{lelievre@cermics.enpc.fr}
}
%
%
\maketitle

\abstract{We present how entropy estimates and logarithmic Sobolev inequalities on the one hand, and the notion of quasi-stationary distribution on the other hand, are useful tools to analyze metastable overdamped Langevin dynamics, in particular to quantify the degree of metastability. We discuss  the interest of these approaches to estimate the efficiency of some classical algorithms used to speed up the sampling, and to evaluate the error introduced by some coarse-graining procedures. This paper is a summary of a plenary talk given by the author at the ENUMATH 2011 conference.}

\section{Introduction and motivation}
\label{sec:intro}

\begin{figure}
\centering\includegraphics[width=5cm]{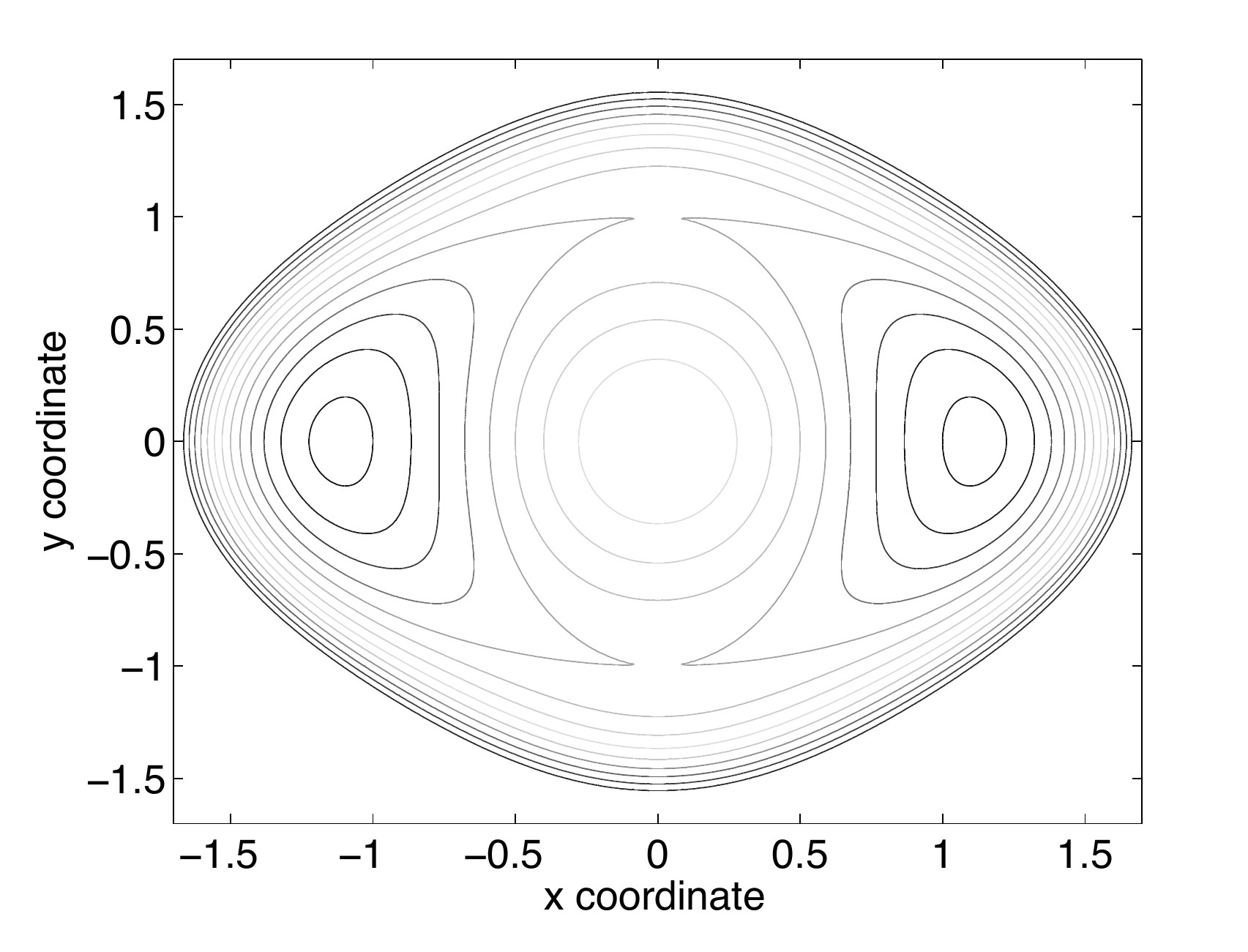}\includegraphics[width=5cm]{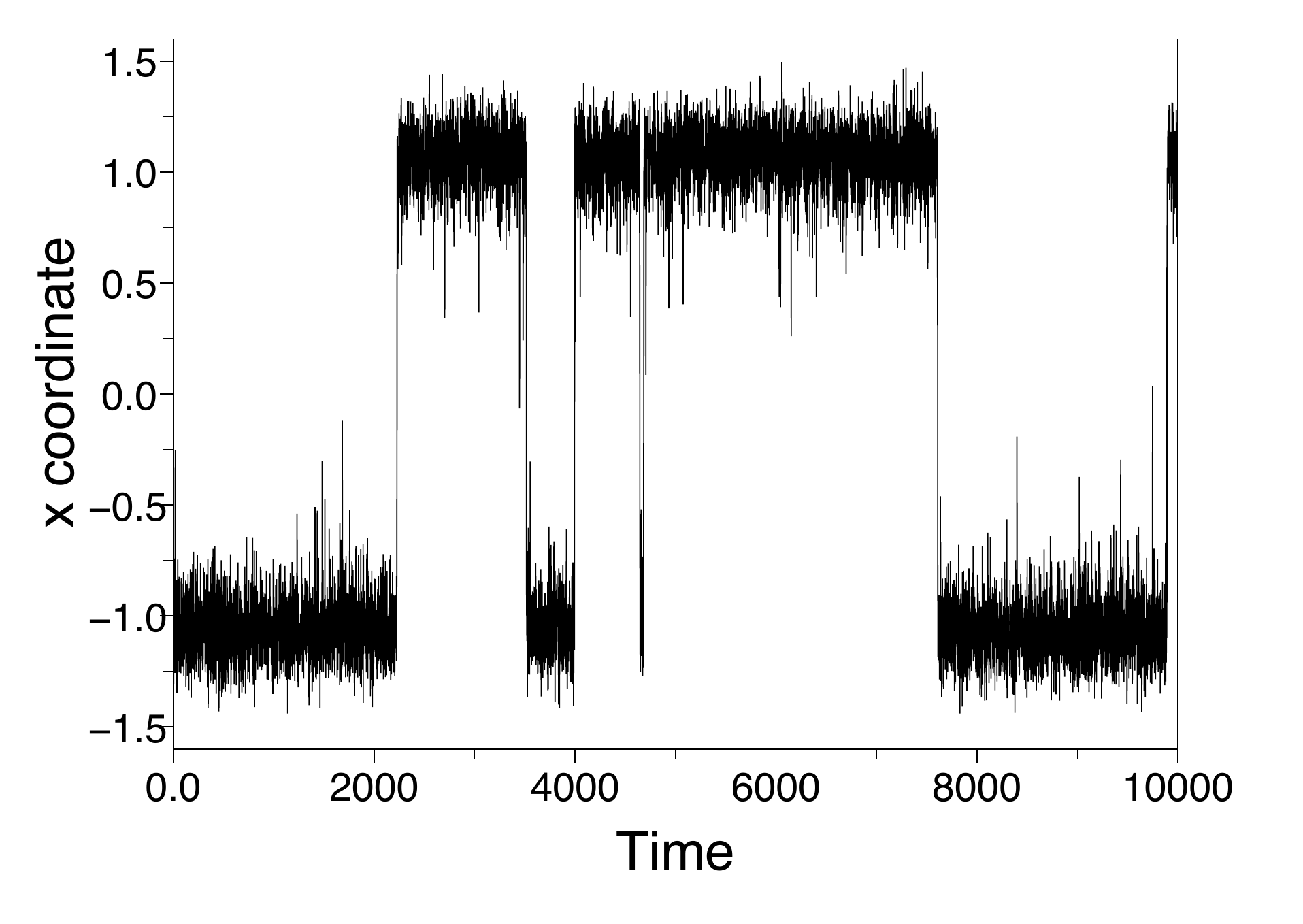}
\centering\includegraphics[width=5cm]{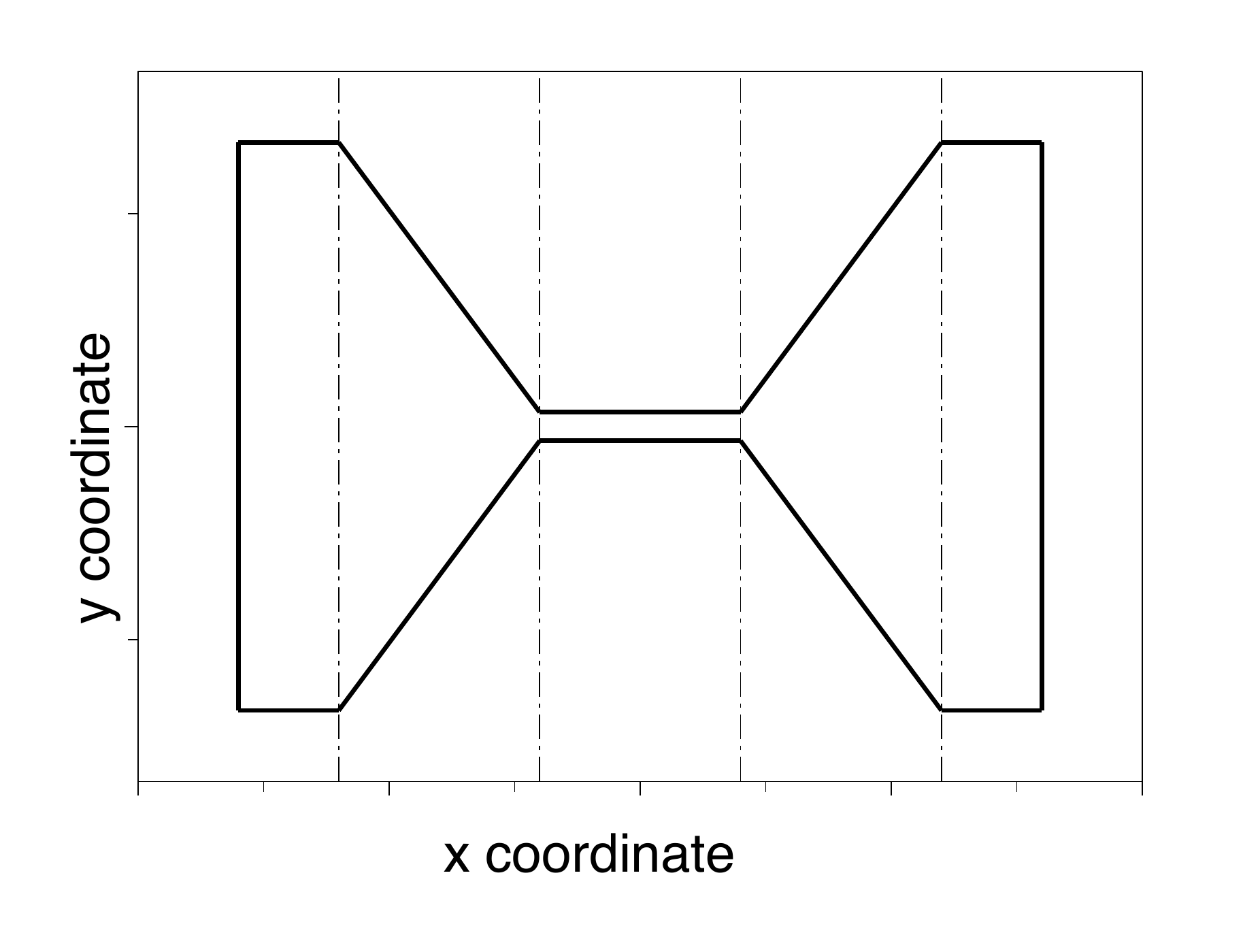}\includegraphics[width=5cm]{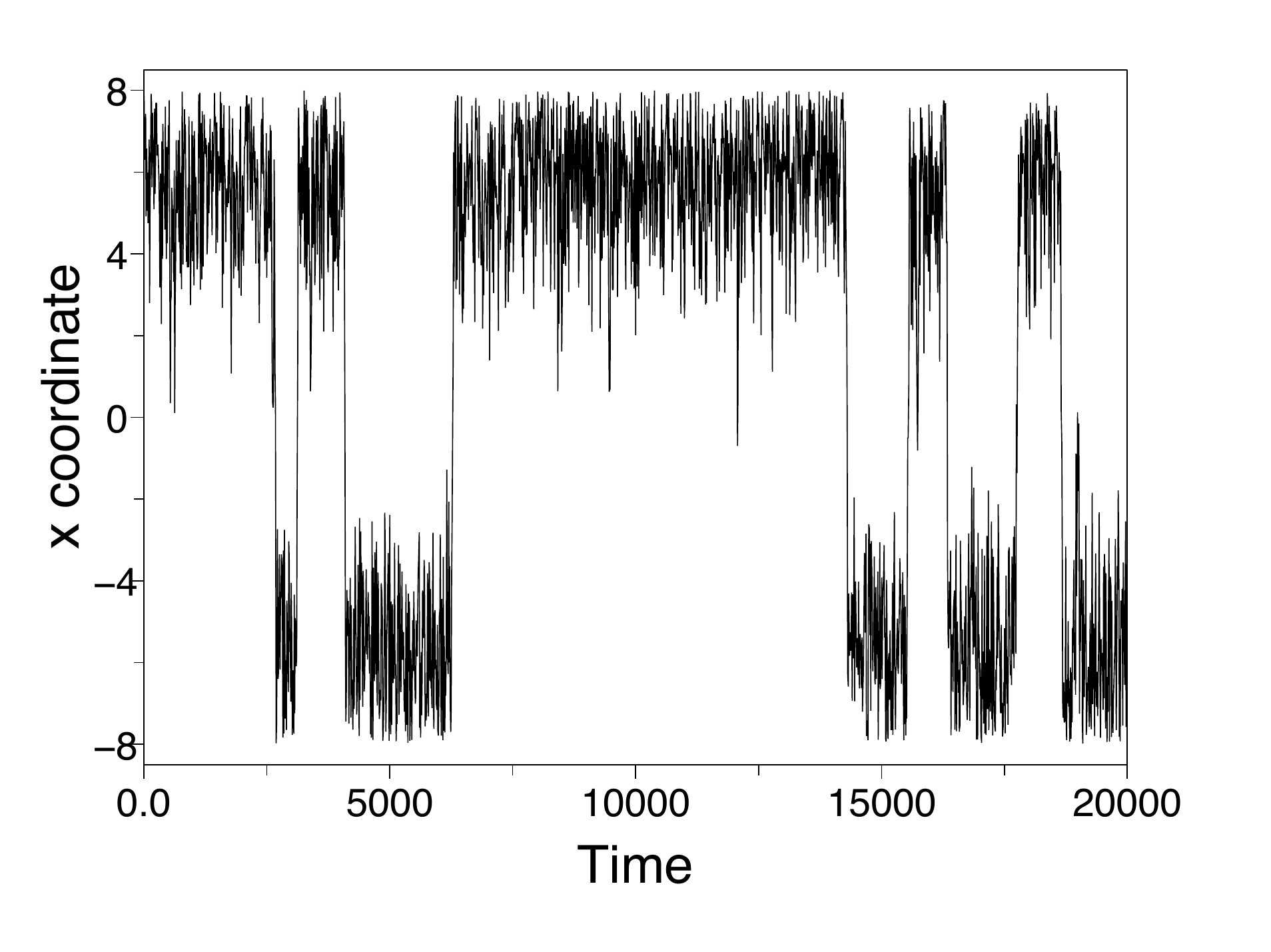}
    \caption{Above: (left) an example of a 2 dimensional potential, with energetic barriers (there are actually two possible saddle points to go from left to right), with (right) the $x$-component as a function of time of the associated stochastic process solution to~\eqref{eq:sde}. Below (left) an example of a 2 dimensional potential (which is zero inside the closed solid-line shape, and infinite outside) for which an entropic barrier (right) is observed.}
    \label{fig:ener_entr}
\end{figure}

The aim of this paper is to present two mathematical viewpoints on metastability. Roughly speaking, a dynamics is said to be metastable if it spends a lot of time in a region (called a metastable state) before hopping to another region. To be more specific, we will focus in the following on the overdamped Langevin dynamics:
\begin{equation}\label{eq:sde}
dX_t = -\nabla V (X_t) \, dt + \sqrt{2 \beta^{-1}} dW_t
\end{equation}
which is used for example in molecular dynamics to describe the evolution of a molecular system. In this context, the configuration of the system $X_t \in \R^n$ is the coordinates of the particles (think of the atoms of a large molecule), $V:\R^n \to \R$ is the potential energy, which to a configuration $x \in \R^n$ associates its energy $V(x)$, and $\beta^{-1}=k_B T$ is proportional to the temperature ($k_B$ being the Boltzmann constant). The stochastic process $W_t$ is a standard $n$-dimensional Brownian motion. For such a dynamics, metastability typically originates from two mechanisms. First, in the small temperature regime, the dynamics~\eqref{eq:sde} can be seen as a perturbation of the simple gradient dynamics $\dot{y}=-\nabla V(y)$ for which, from any initial condition, the solution converges to a local minimum of~$V$. Having this in mind, the dynamics~\eqref{eq:sde} is metastable because it takes a lot of time to leave the vicinity of a local minimum before jumping to the neighborhood of another local minimum. This is due to the {\em energy barriers} which have to be overcome. Such barriers and the zero temperature limit can be analyzed in particular with large deviation techniques~\cite{freidlin-wentzell-84}. Second, metastability may come from entropic effects. Imagine that the configuration space is made of two boxes linked by a narrow corridor (the potential $V$ is  zero on this configuration space, and infinite outside). Then, the dynamics~\eqref{eq:sde} is metastable because it takes a lot of time to find the small corridor to go from one box to the other. Metastability is here due to {\em entropic barriers}. This can be quantified using the notion of free energy, see~\cite{lelievre-rousset-stoltz-book-10,legoll-lelievre-10}. In practice, metastability thus originates from a combination of both energetic and entropic effects, the relative importance of each depending on the system under consideration, and on the temperature. See Figure~\ref{fig:ener_entr} for a numerical illustration. Generally speaking, metastability is related to the multimodality of the measure $\mu$, namely the fact that some high probability regions are separated by low probability regions.

An important concept related to metastability is ergodicity. Under adequate assumptions on $V$, the dynamics~\eqref{eq:sde} can be shown to be ergodic with respect to the canonical measure:
\begin{equation}\label{eq:mu}
\mu(dx)=Z^{-1} \exp(-\beta V(x)) \, dx,
\end{equation}
where $Z = \int_{\R^n} \exp(-\beta V(x)) \, dx$ is assumed to be finite. Ergodicity actually refers to two different properties: (i) an average along a trajectory converges in the long-time limit towards an average with respect to $\mu$: for any test functions $\varphi: \R^n \to \R$,
\begin{equation}\label{eq:erg1}
\lim_{T \to \infty} \frac{1}{T} \int_0^T \varphi(X_t) \, dt  = \int_{\R^n} \varphi d\mu,
\end{equation}
and (ii) the law of the process $X_t$ at time $t$ converges to $\mu$ in the long-time limit: for any test functions $\varphi: \R^n \to \R$,
\begin{equation}\label{eq:erg2}
\lim_{T \to \infty} \E(\varphi(X_T)) = \int_{\R^n} \varphi d\mu.
\end{equation}
If $X_t$ is metastable, both limits~\eqref{eq:erg1} and~\eqref{eq:erg2} are typically very difficult to reach, since $T$ should be sufficiently large to visit all the metastable states. From a numerical viewpoint, metastability raises thus sampling issues, both to compute canonical averages (namely averages with respect to $\mu$) and to compute averages over paths which requires to generate efficiently metastable dynamics, the latter being of course a more complicated task than the former.

In the following, we would like to introduce two mathematical tools to measure the ``degree of metastability'' of the dynamics~\eqref{eq:sde}. In Section~\ref{sec:LSI}, we discuss the notion of {\em Logarithmic Sobolev inequality}, which is a way to quantify the ergodic features of a process, and more precisely, how fast the convergence~\eqref{eq:erg2} happens. In this context, the slower the convergence, the more metastable the process is. In Section~\ref{sec:QSD}, we introduce the notion of {\em quasi-stationary distribution}, and identify the typical time it takes, for a given region of the state space, to reach a quasi equilibrium in this region before leaving it. Metastability in this case is related to the fact that this time is small compared to the typical time it takes to leave the region. In both cases, we will explain how these tools can be used to (i) analyze numerical methods which are used in molecular dynamics to ``accelerate'' metastable dynamics (see Sections~\ref{sec:ABF} and~\ref{sec:ParRep}) and (ii) obtain coarse-grained descriptions of the metastable dynamics (see Sections~\ref{sec:EffDyn} and~\ref{sec:kMC}).

We would like to emphasize the importance of quantifying metastability for practical aspects. Indeed, there exists in the literature many asymptotic analysis in some limiting regimes (zero temperature limit, time-scale separation limit enforced through an explicit small parameter introduced in the dynamics, etc...) where it is shown how a metastable dynamics converges to some effective Markovian dynamics. In practice, the parameters which are considered to go to zero in the asymptotic regimes may indeed be small, but are certainly not zero. A natural question is to quantify the error introduced by assuming that these parameters are zero, an assumption which is behind many numerical methods, and coarse-graining approaches. This requires in turn quantifying the metastable features of the original dynamics. This is precisely the aim of both approaches presented below.

Let us finally clearly state that the aim of this paper is not to provide proofs of the announced results and sometimes even not to state them very precisely mathematically, but to gather in a new and hopefully enlightening way various recent studies, in particular~\cite{lelievre-rousset-stoltz-08,legoll-lelievre-10,le-bris-lelievre-luskin-perez-11}. All the statements below can be reformulated as precise mathematical claims, with rigorous proofs, except the discussion in Section~\ref{sec:kMC} which is more prospective.

\begin{remark}
There are other techniques to quantify metastability, that we do not review here. We would like to mention in particular spectral approaches~\cite{helffer-nier-06,sarich-noe-schuette-10}, potential theoretic approaches~\cite{bovier-eckhoff-gayrard-klein-04,bovier-gayrard-klein-05} and approaches based on drift conditions~\cite{huisinga-meyn-schuette-04}. Drawing connections between these various techniques is an interesting subject, see for example~\cite{cattiaux-guillin-10} for connections between approaches based on drift conditions, and functional inequalities such as LSI.
\end{remark}

\begin{remark}
In this paper, we concentrate on the overdamped Langevin dynamics~\eqref{eq:sde} even though this is not the most widely used dynamics in molecular dynamics. All the algorithms we present below generalize to (and are used with) the phase-space Langevin dynamics, which is much more popular. However, generalizing the mathematical approaches outlined below to the Langevin dynamics is not an easy task, due to the lack of ellipticity of the associated infinitesimal generator, see~\cite{helffer-nier-05,villani-09} for examples of studies in that direction.
\end{remark}

\section{Logarithmic Sobolev inequality}
\label{sec:LSI}

As explained above,  we  quantify in this section the metastability of~\eqref{eq:sde} by considering the rate of convergence of the limit~\eqref{eq:erg2}. We thus consider the law at time~$t$ of~$X_t$, which has a density $\psi(t,x)$. The probability density function $\psi$ satisfies the Fokker-Planck equation:
\begin{equation}\label{eq:FP}
\partial_t \psi = \div( \nabla V \psi + \beta^{-1} \nabla \psi).
\end{equation}
Notice that the density of $\mu$ with respect to the Lebesgue measure, denoted by 
$\psi_\infty(x)=Z^{-1} \exp(-\beta V(x))$, is obviously a stationary solution to~\eqref{eq:FP}.

\subsection{Definition}

Let us introduce the notion of logarithmic Sobolev inequality (see for example~\cite{ABC-00}).
\begin{definition}
The probability measure $\mu$ is said to satisfy a logarithmic Sobolev inequality with constant $R$ (in short LSI($R$)) if and only if, for any probability measure $\nu$ such that $\nu$ is absolutely continuous with respect to $\mu$,
$$H(\nu | \mu) \le \frac{1}{2R} I(\nu | \mu)$$
where $H(\nu|\mu) = \int_{\R^n} \ln\left( \frac{d\nu}{d\mu} \right) \, d\nu$ is the relative entropy of $\nu$ with respect to $\mu$, and $I(\nu|\mu) = \int_{\R^n} \left| \nabla \ln\left( \frac{d\nu}{d\mu} \right) \right|^2 \, d\nu$ is the Fisher information of $\nu$ with respect to $\mu$.
\end{definition}
When both $\mu$ and $\nu$ admit densities (respectively $\psi$ and $\phi$) with respect to the Lebesgue measure, we shall also use the notation $H(\psi|\phi)$ for $H(\mu|\nu)$. A crucial property is the following:
\begin{proposition}\label{prop:cv_lsi}
The measure $\mu$ defined by~\eqref{eq:mu} satisfies a LSI($R$) if and only if, for all probability density functions $\psi_0$, for all time $t \ge 0$,
\begin{equation}\label{eq:exp_CV}
H(\psi(t,\cdot)|\psi_\infty) \le H(\psi_0|\psi_\infty) \exp(- 2 \beta^{-1} R t)
\end{equation}
where $\psi$ is the solution to~\eqref{eq:FP} with initial condition $\psi(0,\cdot)=\psi_0$.
\end{proposition}
This is a simple consequence of the standard computation: if $\psi$ satisfies~\eqref{eq:FP}, then
$$\frac{d}{dt} H(\psi(t,\cdot)|\psi_\infty)=-\beta^{-1}I(\psi(t,\cdot)|\psi_\infty).$$
A natural way to quantify metastability if thus to relate it to $R$:
\begin{equation}\label{eq:metastab_LSI1}
\text{The smaller $R$, the more metastable the dynamics~\eqref{eq:sde} is.}
\end{equation}

Before we proceed, let us make two remarks. First, by combining the classical Bakry-Emery criteria with the perturbation result of Holley and Stroock, the measure $\mu$ actually satisfies a LSI under very mild assumptions on $V$: basically, if $V$ is smooth, and is $\alpha$-convex at infinity, then a LSI for $\mu$ holds. What is more complicated is to get the optimal constant $R$. Second, in the simple case of a double well potential in dimension 1, it is standard to show that the average time spent by the process $X_t$ in a well before hopping to another one increases as the temperature goes to zero, using for example large deviation techniques. It can also be checked in this simple setting that the constant $R$ scales like $\exp(- \beta H)$ in the limit of large $\beta$ (small temperature), where $H$ is the height of the barrier to overcome to leave a given well. This prototypical situation thus shows that the LSI constant indeed allows to quantifying the  intuitive definition of metastability we gave in the introduction.

\subsection{Metastability along a reaction coordinate}

Many algorithms and modelling discussions are based on the introduction of a so-called reaction coordinate $\xi$, namely a smooth low-dimensional function which typically indices transition from one metastable state to another. For simplicity, let us assume that $\xi:\R^n \to \T$ has values in the one-dimensional torus (think of an angle in a molecule, which characterizes its conformation). Then, one may introduce probability measures associated to $\mu$ and $\xi$:
\begin{itemize}
\item the image $\xi*\mu$ of the measure $\mu$ by $\xi$, which is a probability measure on the torus $\T$, and which is also written as $\xi*\mu(dz) = \exp(-\beta F(z)) \, dz$, $F: \T \to \R$ being then the so-called free energy associated to $\mu$ and $\xi$. Using the co-area formula, a formula for $F$ is given by
\begin{equation}\label{eq:F}
F(z) = -\beta^{-1} \ln \int_{\Sigma(z)} Z^{-1} \exp(-\beta V(x)) \delta_{\xi(x) -z}(dx)
\end{equation}
where $\Sigma(z)=\{ x \in \R^n, \xi(x) = z\}$ and $\delta_{\xi(x)-z}(dx)$ is a measure supported by $\Sigma(z)$ such that $dx = \delta_{\xi(x)-z}(dx) \, dz$.
\item the family of conditional probability measures $\mu(\cdot | \xi(x)=z)$ with support $\Sigma(z)$, which are indexed by $z \in \T$ and defined as
$$\mu(dx | \xi(x)=z) = \frac{\exp(-\beta V(x)) \delta_{\xi(x) -z}(dx)}{\exp(-\beta F(z))}.$$
\end{itemize}
These two measures are completely defined through the conditioning formula: for any test functions: $\varphi: \R^n \to \R$ and $\psi: \T \to \R$,
$$\int_{\R^n} \varphi(\xi(x)) \psi(x) \mu(dx) = \int_{\T} \varphi(z) \int_{\Sigma(z)} \psi(x) \mu(dx| \xi(x)=z) \, \exp(-\beta F(z)) \, dz.$$
For proofs, we refer for example to~\cite{lelievre-rousset-stoltz-book-10}.

Let us assume that the measures $\mu(\cdot | \xi(x)=z)$ satisfy a LSI($\rho$) (with $\rho$ uniform in $z$), and that the measure $\xi*\mu=Z^{-1} \exp(-\beta F(z)) \, dz$ satifies a LSI($r$). In the spirit of the previous subsection, we will say that (this will be assumption (H1) below):
\begin{equation}\label{eq:metastab_LSI2}
\text{``the metastability of the process $X_t$ is along $\xi$'' if and only if $\rho \gg r$.}
\end{equation}
A typical example of such a situation is given on Figure~\ref{fig:ener_entr}, for which $\xi(x,y)=x$. In such a two-dimensional setting, notice that $F(x)=-\beta^{-1} \ln \int_{\R} Z^{-1} \exp(-\beta V(x,y)) \, dy$ and $\displaystyle{\mu(dy | \xi(x,y)=x) = \frac{\exp(-\beta V(x,y)) \, dy}{\exp(-\beta F(x))}}$.

It is possible to relate the LSI constant of the measure $\mu$ (namely $R$) to the LSI constants of the measures  $\mu(\cdot | \xi(x)=z)$ (namely $\rho$) and $\xi*\mu$ (namely $r$), see~\cite{lelievre-09}. Roughly speaking, 
\begin{equation}\label{eq:JFA}
\text{If $r$ is small and $\rho$ is small, then $R$ is small.}
\end{equation}
Moreover, if $R$ is small and $\xi$ is well chosen, then $r$ is small but $\rho$ may be very large. In such a case, the metastability of the process $X_t$ is essentially encoded in the low-dimensional observable $\xi(X_t)$. It is then possible to use numerical and analytical techniques to accelerate the long-time convergence and go around the difficulties associated to metastability, as explained in the two next subsections. This could also be used to yield a definition of what a good reaction coordinate is: it is a (low-dimensional) function $\xi$ such that $\rho/r$ is as large as possible. Designing a numerical method which would look for the best $\xi$ in this respect would be very interesting for practical applications.

Before we proceed, we provide a formula that will be useful below. By using the co-area formula, starting from~\eqref{eq:F}, it is possible to check that the derivative of $F$ writes:
\begin{equation}\label{eq:F'}
F'(z)=\int f(x) \mu(dx | \xi(x)=z) \text{ with } f=\frac{\nabla V\cdot \nabla \xi}{|\nabla \xi|^2} - \beta^{-1} \div\left(\frac{\nabla \xi}{|\nabla \xi|^2} \right).
\end{equation}
In the simple case $\xi(x,y)=x$ mentioned above, the function $f$ is simply $f=\partial_x V$.

\subsection{A first example: the adaptive biasing force technique}
\label{sec:ABF}

As explained above, one difficulty related to metastability is that the convergence~\eqref{eq:erg2} is very slow. In particular, this implies that it will be difficult to sample the canonical measure $\mu$ from a trajectory $X_t$. The fact that it is difficult to sample a multimodal measure is a well known problem shared with other fields than molecular dynamics. In statistics for example, Markov Chain Monte Carlo methods are very popular and similar sampling problems occur for Bayesian inference~\cite{chopin-lelievre-stoltz-11}.

\begin{figure}
\centering\includegraphics[width=5cm]{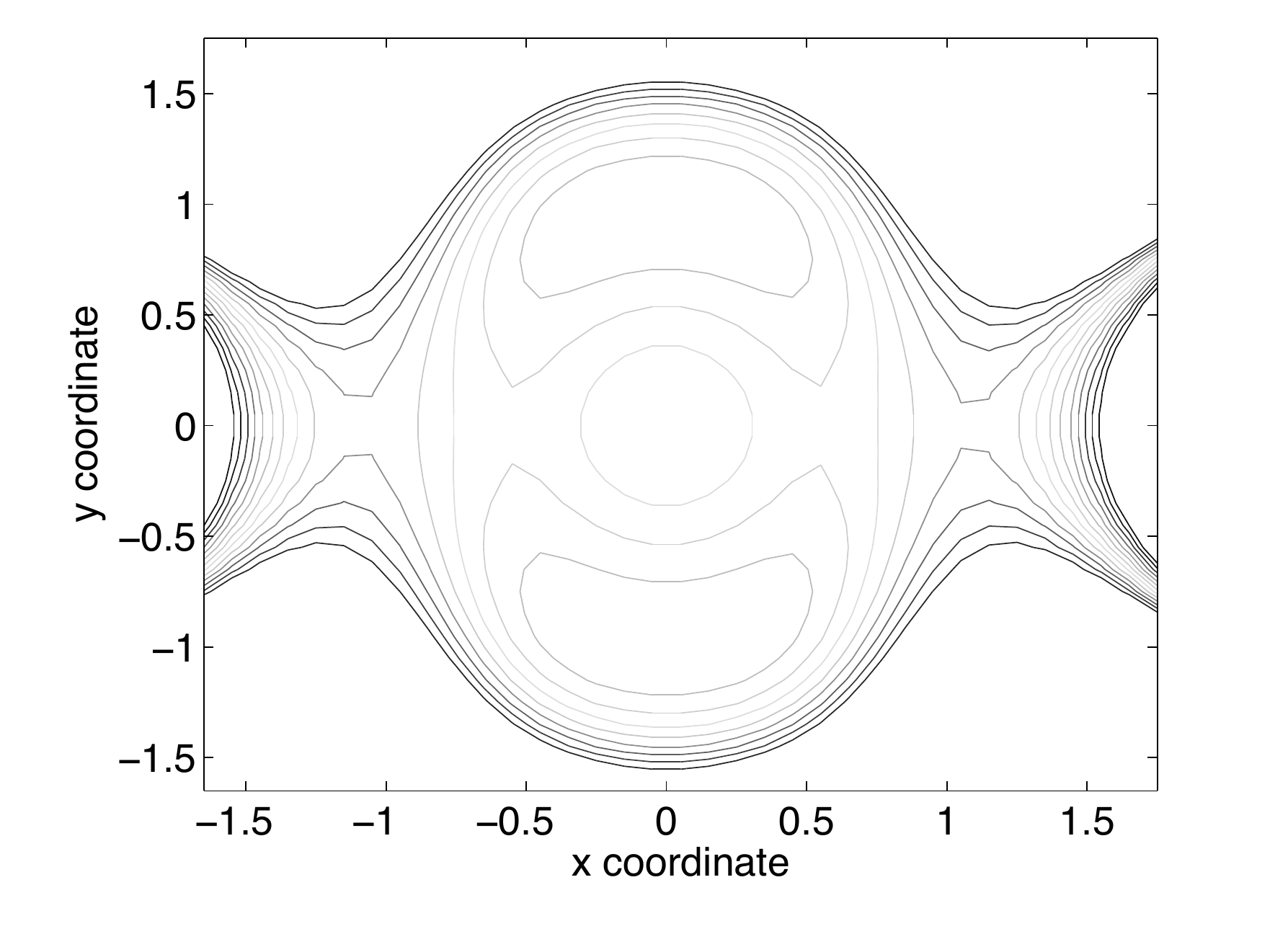}\includegraphics[width=5cm]{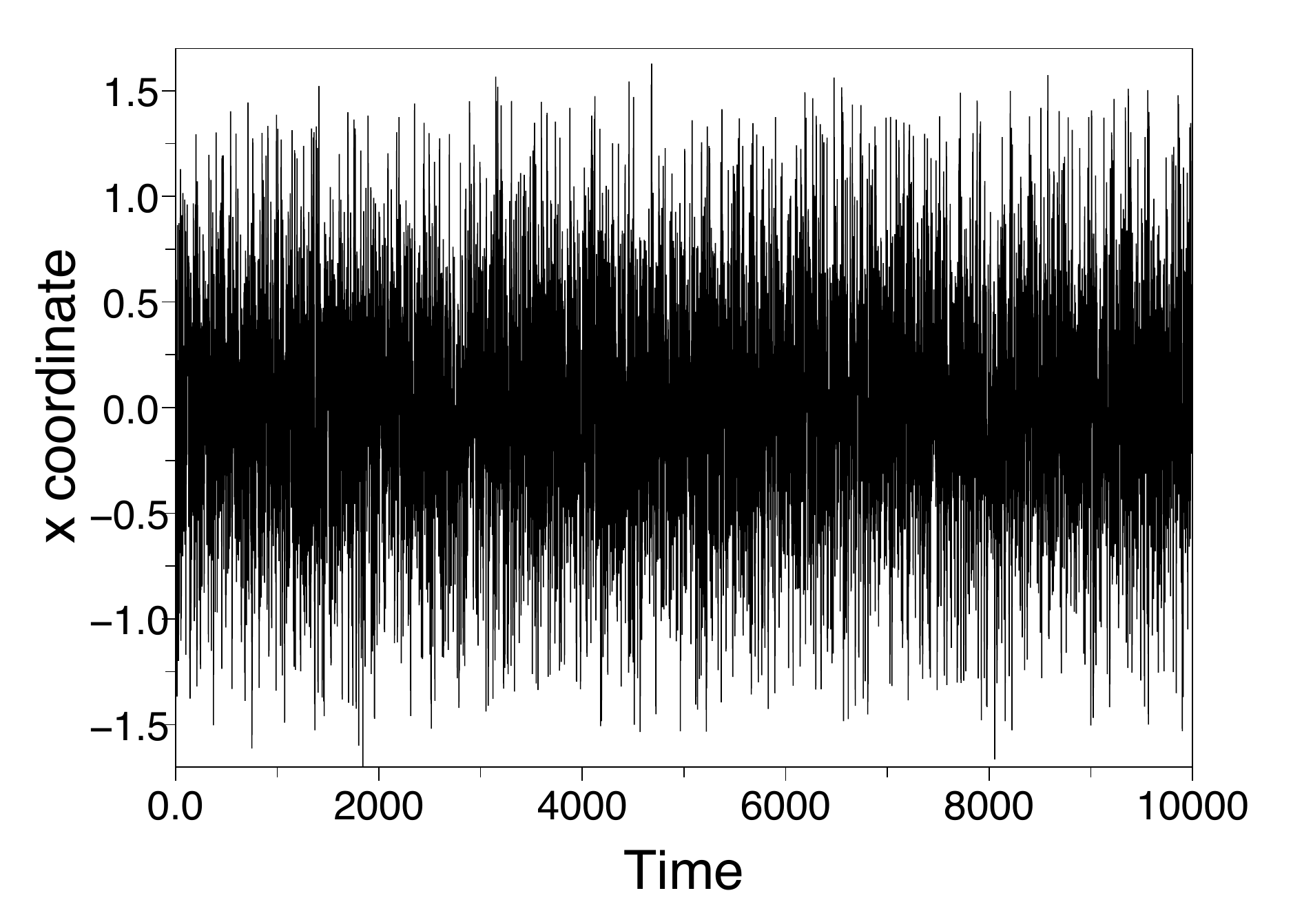}
\centering\includegraphics[width=5cm]{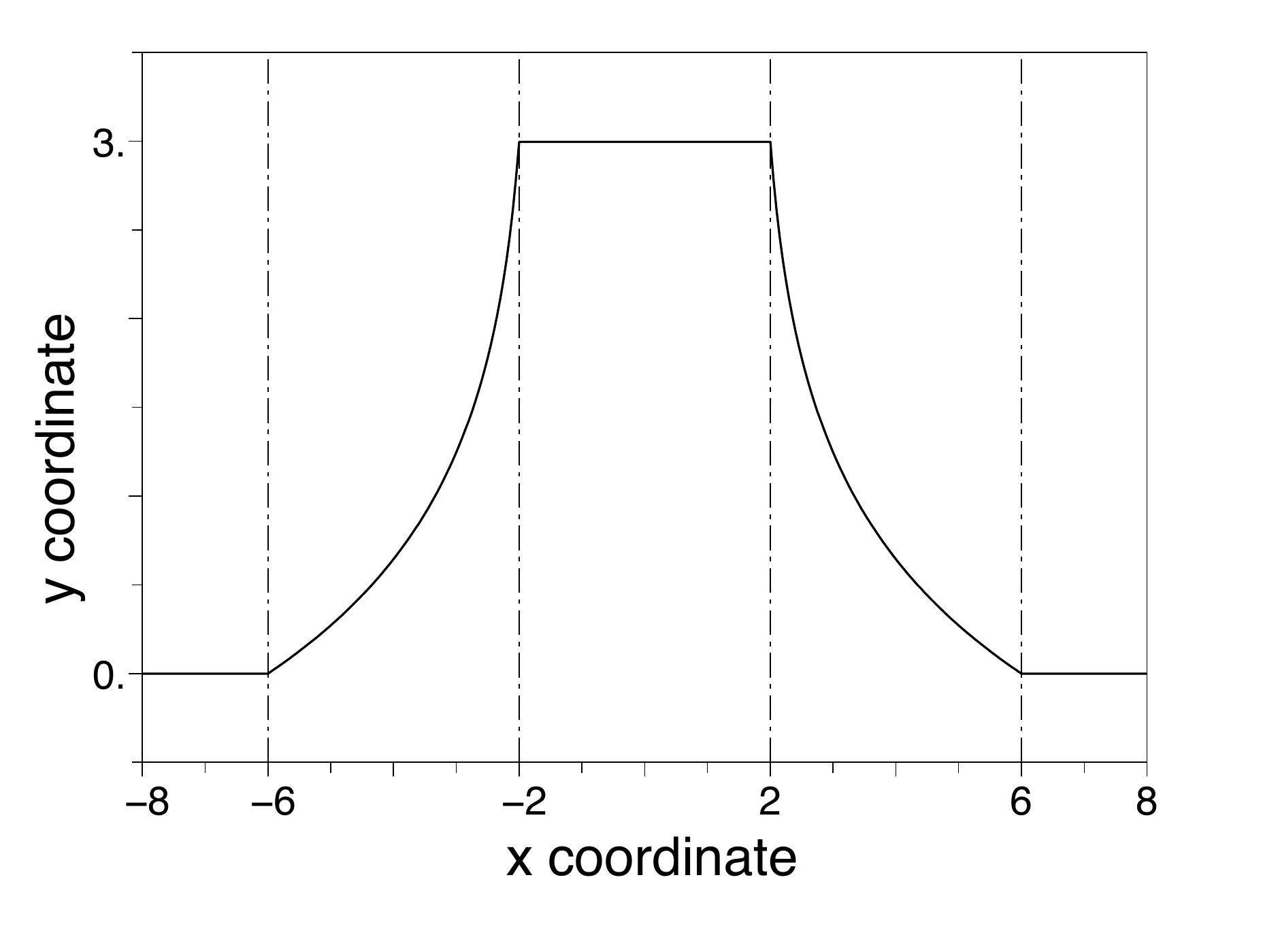}\includegraphics[width=5cm]{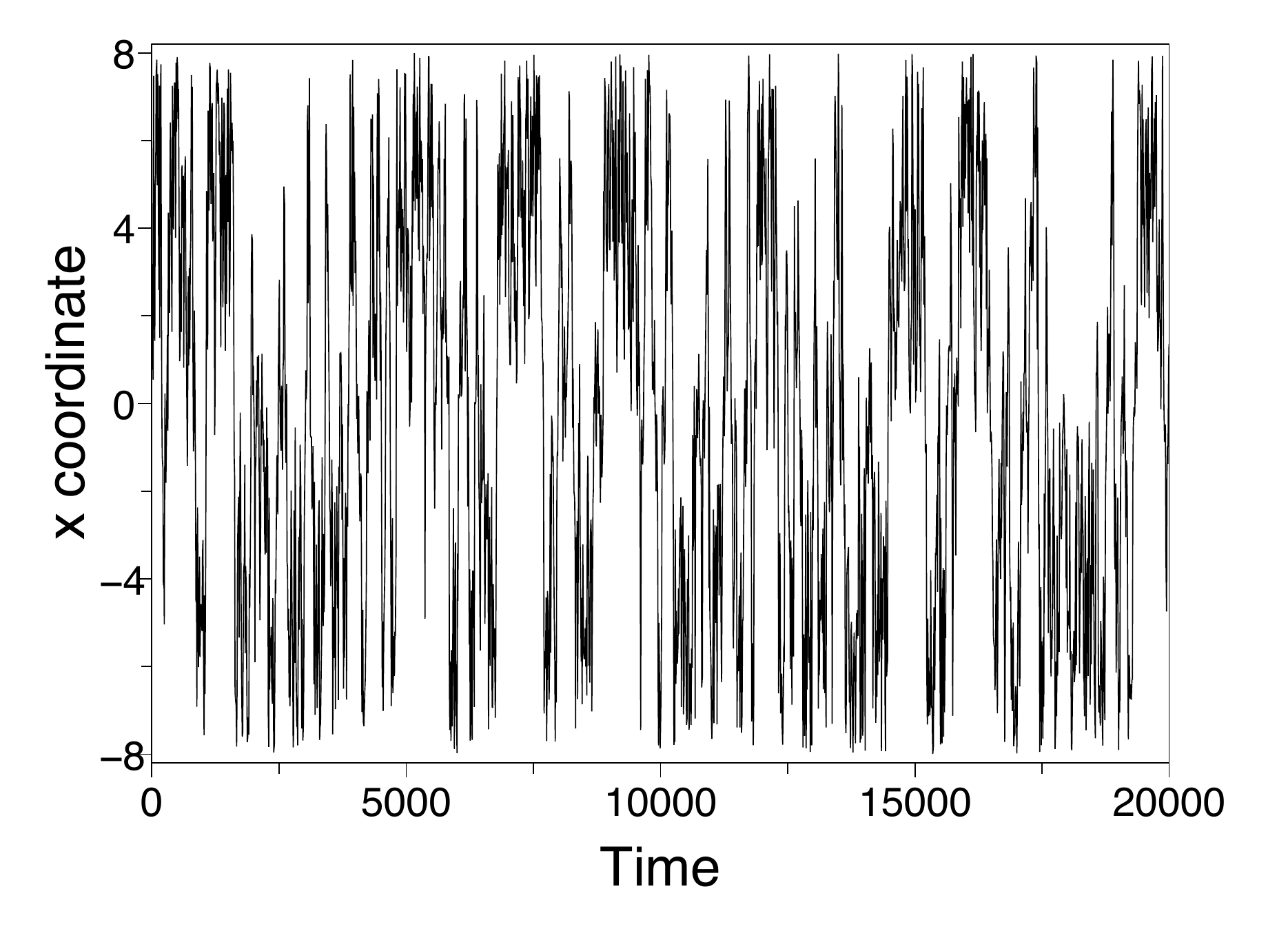}
    \caption{See Figure~\ref{fig:ener_entr} for comparison. The reaction coordinate is $\xi(x,y)=x$. Above: In the energetic barrier case, (left) the 2 dimensional potential minus the free energy and (right) the $x$-component of the associated stochastic process. Below: In the entropic barrier case, (left) the free energy and (right) the $x$-component of the stochastic process simulated again with the free energy biased potential.}
    \label{fig:FE_bias}
\end{figure}

Under assumption~\eqref{eq:metastab_LSI2}, a natural importance sampling idea is the following. If the metastability of $X_t$ is along $\xi$, it is sensible to try to remove these metastable features by changing the potential $V$ to $V-F\circ \xi$, where $F$ is the free energy~\eqref{eq:F} (and $F \circ \xi$ denotes the composition of $F$ with $\xi$). Indeed, if we denote $$\mu_F(dx)=Z_F^{-1} \exp(-\beta (V - F\circ \xi) (x)) \, dx$$ the associated tilted measure, we clearly have $\xi * \mu_F = 1_{\T}(x) \, dx$ and $\mu_F(\cdot | \xi(x)=z)=\mu(\cdot | \xi(x)=z)$: the conditional measures remain the same, but the marginal along $\xi$ is now a uniform measure, namely a very gentle measure, without any multimodality, and thus easy to sample. In other words, following~\eqref{eq:JFA}, we may hope that the LSI constant of the tilted measure $\mu_F$ is much smaller than the one of the original measure~$\mu$. As a numerical illustration of this fact, we present on Figure~\ref{fig:FE_bias} the equivalent of the trajectories presented on Figure~\ref{fig:ener_entr}, using the biased potential $V-F\circ \xi$. We clearly observe much more transitions from left to right with the biased potential, both in the case of an energetic barrier and an entropic barrier.

The difficulty is of course that computing the free energy $F$ is a challenge in itself (it corresponds to a sampling problem of, {\em a priori}, a similar complexity as the sampling of $\mu$). The free energy is actually a quantity of great interest for practitioners in itself~\cite{chipot-pohorille-07}, and many methods have been designed to compute it~\cite{lelievre-rousset-stoltz-book-10}. Thus, biasing the measure using $F$ does not seem to be such a good idea. The principle of adaptive biasing technique is to actually use, at a given time $t$, {\em an approximation} $F_t$ of the free energy $F$ (in view of the configurations visited so far) in order to bias the dynamics. In other words, instead of using the biased dynamics (which assumes that $F$ is already known)
\begin{equation}\label{eq:FE_bias_sde}
\left\{
\begin{aligned}
dX_t&=-\nabla (V- F \circ \xi) (X_t) \, dt + \sqrt{2 \beta^{-1}} d W_t,\\
F'(z)&=\E_\mu(f(X) | \xi(X)=z),
\end{aligned}
\right.
\end{equation}
where the formula for $F'(z)$ is exactly~\eqref{eq:F'} (here and in the following, $\E_\mu$ denoting an expectation taken with respect to the measure $\mu$), one rather considers
\begin{equation}\label{eq:ABF}
\left\{
\begin{aligned}
dX_t&=-\nabla (V- F_t \circ \xi) (X_t) \, dt + \sqrt{2 \beta^{-1}} d W_t,\\
F'_t(z)&=\E(f(X_t) | \xi(X_t)=z).
\end{aligned}
\right.
\end{equation}
The bottom line is that if $X_t$ solution to~\eqref{eq:ABF} was at equilibrium instantaneously with respect to $\mu_{F_t}$, then $F_t'$ would be exactly $F'$. Of course, this instantaneous equilibrium assumption is wrong, but the hope is that $F_t'$ learns on the fly (and eventually converges to) $F'$.

The dynamics~\eqref{eq:ABF} is the so called adaptive biasing force (ABF) process, and is one of the most efficient methods to compute free energy differences, see~\cite{darve-pohorille-01,henin-chipot-04} for the original idea. One way to understand such a method is that already visited states are penalized (the potential is flattened in these regions, or equivalently, the probability to visit them is increased) in order to force the stochastic process to visit new regions. There are many other techniques along these lines, see~\cite{lelievre-rousset-stoltz-07-b}. The hope is that, by forcing the system to visit all the possible values of $\xi$, the metastability of the original dynamics is completely overcome. This is somewhat in the spirit of simulated annealing or parallel tempering, where the temperature (which in some sense plays the role of $\xi$) of the system is changed in order to visit new regions.

What can be shown is that under the assumptions:
\begin{itemize}
\item (H1) the metastability of the process $X_t$ is along $\xi$ (see~\eqref{eq:metastab_LSI2}),
\item (H2) the cross derivative $\nabla_{\Sigma(z)} f$ is bounded (where $\nabla_{\Sigma(z)}$ denotes the gradient projected onto the tangent space to $\Sigma(z)$),
\end{itemize}
then the convergence of the ABF process~\eqref{eq:FE_bias_sde} to equilibrium is much faster (basically exponential with rate $\beta^{-1} \rho$) compared to the convergence of the original process~\eqref{eq:sde} to equilibrium (which is exponential with rate $\beta^{-1} R$, see~\eqref{eq:exp_CV}). In particular, $F_t$ converges to the free energy $F$ very fast if $\xi$ is such that $\rho$ is large (this is assumption (H1)).
We refer to~\cite{lelievre-rousset-stoltz-08} for a precise mathematical statement. The proof is based on entropy techniques~\cite{arnold-markowich-toscani-unterreiter-01}, and the idea of two-scale analysis for logarithmic Sobolev inequalities~\cite{grunewald-otto-villani-westdickenberg-09,lelievre-09}.

We also refer to~\cite{lelievre-minoukadeh-11} for some refinements in the cases when $\rho$ is only large for some of the family of conditional probability measures $\mu(\cdot | \xi(x)=z)$ indexed by $z$ (the so-called bi-channel situation). For practical aspects (discretization techniques and numerical illustrations), we refer to~\cite{lelievre-rousset-stoltz-07-b,chipot-lelievre-11,jourdain-lelievre-roux-10} and to~\cite{chopin-lelievre-stoltz-11} for an application of such techniques in the context of Bayesian inference in statistics.

\begin{remark}
In the long-time limit, $F$ is obtained, and the measure sampled by~\eqref{eq:ABF} is not $\mu$ but $\mu_F$. There are basically two ways to recover averages with respect to $\mu$. First, as in standard importance sampling approaches, reweighting can be used:
$$\E_{\mu} (\varphi(X))=\frac{\E_{\mu_F} (\varphi(X) \exp(-\beta F \circ \xi (X))}{\E_{\mu_F} (\exp(-\beta F \circ \xi (X))}.$$
For this idea to be efficient, the weights should not be too widespread (otherwise the variance may be large), which means that $\sup F - \inf F$ should not be too large, see~\cite{chopin-lelievre-stoltz-11} for a discussion of these aspects.
Another idea is to use a conditioning approach:
$$\E_{\mu} (\varphi(X))=\frac{\displaystyle \int_{\T} \E_\mu(\varphi(X)|\xi(X)=z) \exp(-\beta F(z)) \, dz}{\displaystyle\int_{\T}  \exp(-\beta F(z)) \, dz}.$$
The conditional probabilities $\E_\mu(\varphi(X)|\xi(X)=z)$ can then be computed using either the fact that $\E_\mu(\varphi(X)|\xi(X)=z)=\E_{\mu_F}(\varphi(X)|\xi(X)=z)$ (so that the ABF process~\eqref{eq:ABF} can be used to compute them) or using dedicated techniques to sample the conditional probability measure $\mu(\cdot|\xi(x)=z)$, which should be easy under assumption (H1). Indeed, if $\rho$ is large, in virtue of Proposition~\ref{prop:cv_lsi}, the overdamped Langevin dynamics associated to the measure~$\mu(\cdot|\xi(x)=z)$ should converge very fast to equilibrium. Such a dynamics is roughly speaking a projection of the original gradient dynamics~\eqref{eq:sde} to the submanifold~$\Sigma(z)$. We refer to~\cite{ciccotti-lelievre-vanden-einjden-08,lelievre-rousset-stoltz-book-10,lelievre-rousset-stoltz-11} for more information on such constrained sampling techniques.
\end{remark}

\subsection{A second example: obtaining an effective dynamics on $\xi(X_t)$}
\label{sec:EffDyn}

If we are in the situation~\eqref{eq:metastab_LSI2} where the metastability of $X_t$ is along $\xi$, another idea is to try to derive an effective dynamics for $\xi(X_t)$, which would then be easy to simulate since it is low-dimensional, and hopefully associated with a smaller characteristic timescale than the original dynamics~\eqref{eq:sde}. The idea is that if the metastability is along $\xi$, then $\xi(X_t)$ should move much more slowly than the components of $X_t$ which are along ``directions orthogonal to $\xi$'', so that some averaging should be possible along those directions. This is very much in the spirit of projection operator or Mori-Zwanzig techniques~\cite{givon-kupferman-stuart-04}. Below, we first derive an effective Markovian dynamics for $\xi(X_t)$ and then assess the quality of this effective dynamics by deriving quantitative bounds. These quantitative bounds are again obtained using logarithmic Sobolev inequalities and entropy computations.

The idea to derive an effective dynamics on $\xi(X_t)$ starts from a simple It\^o calculus. If $(X_t)_{t \ge 0}$ satisfies~\eqref{eq:sde}, then
$$d\xi(X_t) = ( -\nabla V \cdot \nabla \xi + \beta^{-1} \Delta \xi) (X_t) \, dt + \sqrt{2 \beta^{-1}} |\nabla \xi (X_t)| \frac{\nabla \xi (X_t)}{|\nabla \xi (X_t)|}\cdot dW_t.$$
Of course, this is not a closed equation for the evolution of $\xi(X_t)$.
It is not difficult to check that if we the consider
$$d\tilde{z}_t=\tilde{b}(t,\tilde{z}_t) \, dt + \sqrt{2 \beta^{-1}} \tilde{\sigma} (t,\tilde{z}_t) \, dB_t$$
with
$$\tilde{b}(t,\tilde{z})=\E\left( (-\nabla V \cdot \nabla \xi + \beta^{-1} \Delta \xi)(X_t) \Big| \xi(X_t)=\tilde{z}  \right)$$
and
$$\tilde{\sigma}^2(t,\tilde{z})=\E\left( |\nabla \xi|^2(X_t) \Big| \xi(X_t)=\tilde{z}  \right),$$
then, for all time $t \ge 0$, the law of the random variable $\xi(X_t)$ is {\em equal} to the law of the random variable $\tilde{z}_t$. The difficulty of course is that $\tilde{b}$ and $\tilde{\sigma}$ are intractable numerically, since they are functions depending on time. A natural idea (following the intuition given at the beginning of this section) is that one could replace the conditional expectations defining $\tilde{b}$ and $\tilde{\sigma}^2$ by conditional expectations {\em at equilibrium}, namely:
\begin{equation}\label{eq:z}
dz_t=b(z_t) \, dt + \sqrt{2 \beta^{-1}} \sigma (z_t) \, dB_t
\end{equation}
with
$$b(z)=\E_{\mu}\left( (-\nabla V \cdot \nabla \xi + \beta^{-1} \Delta \xi)(X) \Big| \xi(X)=z  \right)$$
and
$$\sigma^2(z)=\E_{\mu}\left( |\nabla \xi|^2(X) \Big| \xi(X)=z  \right).$$
For related approaches, we refer to~\cite{e-vanden-eijnden-04,maragliano-fischer-vanden-einjden-ciccotti-06,pavliotis-stuart-07}.

Now that we have derived a Markovian dynamics~\eqref{eq:z} on $z_t$, which should be such that $(z_t)_{t \ge 0}$ is close to $(\xi(X_t))_{t \ge 0}$, a natural question is whether we can give some error estimate of some ``distance between the two processes''. What we have shown in~\cite{legoll-lelievre-10} is that under the same assumptions (H1) and (H2) needed for the analysis of the longtime convergence of the ABF dynamics, the relative entropy of the law at time $t$ of $\xi(X_t)$ with respect to the law at time $t$ of $z_t$ is bounded from above {\em uniformly in time} by a constant divided by $\rho^2$. Thus, the larger $\rho$ (this is exactly (H1)), the closer these two probability measures are, for all times. We refer to~\cite{legoll-lelievre-10,legoll-lelievre-12} for a precise mathematical statement. The proof uses basically the same ingredients as for the analysis of the long-time convergence of the ABF process.

\begin{remark}\label{rem:olla}
The result we mention above only concerns the closeness of the {\em marginals in time} of the two processes, namely the laws, at times $t \ge 0$, of $\xi(X_t)$ and $z_t$. Of course, the stochastic processes $(\xi(X_t))_{t \ge 0}$ and $(z_t)_{t \ge 0}$ contain much more information than their marginal in times (think of time correlations, first time of escape from a well, typical paths to go from one well to another, etc...), which are of interest in the context of molecular dynamics simulation. A natural question is thus whether one can prove similar results on the {\em law of the paths}. This can indeed be done, under very close assumptions to (H1) and (H2), on finite time intervals. We refer to~\cite{legoll-lelievre-olla-12}.
\end{remark}

\section{Quasi-stationary distribution}
\label{sec:QSD}

In this second part, we would like to introduce another tool to study metastability, and to show in particular how this tool may be useful to analyze the parallel replica dynamics~\cite{voter-98} which is a numerical method to efficiently generate {\em a trajectory} of a metastable dynamics.

Let us again consider the dynamics~\eqref{eq:sde}, and let us assume that we are given some partition of the state space, for example through an application
$${\mathcal S}: \R^n \to \N$$
which to a given configuration $x$ associates the number ${\mathcal S}(x)$ of the state in which~$x$ is. In some sense, ${\mathcal S}$ could be thought as an equivalent to the reaction coordinate $\xi$ of the previous sections, but with discrete values. In the following, one should think of the states (for $k \in \N$)
$$W_k=\{x \in \R^n, \, {\mathcal S}(x)=k \}$$
as the metastable regions we mentioned in the introduction (think for example of the basins of attraction of the local minima of $V$). Typically, when $X_t$ enters a state $W_k$, we would like it to stay in the state for a long-time before it leaves the state. To formalize this idea, and to quantify the error introduced when this assumption is done in some algorithms, we rely on the notion of quasi-stationary distribution. In the following, we assume that the states $W_k$ are smooth bounded connected subsets of $\R^n$.

\subsection{Definition and two basic properties of the quasi-stationary distribution}

In this section, we consider a given state $W_{k}$ for a fixed index $k$, and we denote by $W$ the interior of $W_k$. Let us introduce the notion of the quasi-stationary distribution associated to the well $W$. We refer to~\cite{cattiaux-collet-lambert-martinez-meleard-san-martin-09,martinez-san-martin-04,mandl-61,collet-martinez-san-martin-95,steinsaltz-evans-07,pinsky-85,ferrari-kesten-martinez-picco-95,ferrari-maric-07,ferrari-martinez-san-martin-96} for general introductions to the quasi-stationary distribution.

\begin{definition}
The quasi-stationary distribution (QSD) associated to the dynamics~\eqref{eq:sde} and the state $W$ is defined as a measure $\nu$ with support in $W$ and such that  $\forall t >0$, $\forall A \subset W,$
$${\nu (A) = \frac{\displaystyle \int_W \P(X_t^{x} \in A , \, t < T_W^{x}) \, \nu(dx) }{\displaystyle \int_W \P( t < T_W^{x}) \, \nu(dx)}},$$
where $X_t^x$ denotes the solution to~\eqref{eq:sde} such that $X_0=x$ and $T_W^x=\inf\{ t \ge 0, X_t^x \not\in W\}$.
\end{definition}
In other words, if $X_0$ is distributed according to $\nu$ and if $(X_s)_{s \ge 0}$ solution to~\eqref{eq:sde} has not left the state $W$ on the interval $[0,t]$, then $X_t$ is also distributed according to $\nu$.

The QSD enjoys two other properties which are very important in practice, given in the two following propositions. To present these two properties, we need to state an intermediate result. Let $L=-\nabla V\cdot \nabla + \beta^{-1} \Delta$ be the infinitesimal generator associated to~\eqref{eq:sde}. Then the density $u_1$ of $\nu$ with respect to $\mu$ (namely such that $\nu(dx) = u_1(x) \mu(dx)$) is the first eigenfunction of $L$ with zero Dirichlet boundary condition on $\partial W$ (and with normalization $\int_W u_1 d\mu=1$). In other words:
\begin{equation}
u_1 \in \argmin_{u \in H^1_{\mu,0}(W)} \frac{\displaystyle \int_W |\nabla u|^2 \, d\mu}{\displaystyle \int_W u^2 \, d\mu}
\end{equation}
and
\begin{equation}
\left\{
\begin{aligned}
Lu_1&=-\lambda_1 u_1 \text{ in } W,\\
u_1&=0\text{ on } \partial W.
\end{aligned}
\right.
\end{equation}
Here $H^1_{\mu,0}(W)=\{ u:W \to \R, \, \int_W \left(|\nabla u|^2 + u^2 \right) \, d\mu< \infty \text{ and } u=0 \text{ on } \partial W\}$. The operator $L$ can be shown to be negative, self-adjoint on $L^2_\mu$, and with a discrete spectrum $(-\lambda_1,-\lambda_2, \ldots, -\lambda_n, \ldots)$ (the eigenvalues are in decreasing order and counted with multiplicity) and associated eigenfunctions $(u_1,u_2, \ldots, u_n, \ldots)$ in $H^1_{\mu,0}(W)$. The first eigenstate can be shown to be non-degenerate ($\lambda_2 > \lambda_1 > 0$), and the first eigenfunction does not vanish on $W$, and can thus be assumed to be positive ($u_1 > 0$ on $W$). Using this spectral decomposition of the infinitesimal generator, and this characterization of the QSD, one can show the following results.
\begin{proposition}\label{prop:longtime_QSD}
Let $X_0$ be distributed according to a distribution with support in $W$, and let us consider $(X_t)_{t \ge 0}$ solution to~\eqref{eq:sde}. Then, the law of $X_t$ conditionally to the fact that the process remained in $W$ up to time $t$ converges to the QSD $\nu$: for any test function $\varphi:W \to \R$,
$$\lim_{t \to \infty} \E(\varphi(X_t) \, | \, t< T_W) = \int_W \varphi d\nu.$$
Moreover, the rate of convergence is exponential:
\begin{equation}\label{eq:cv_qsd}
\sup_{\varphi \in L^\infty(W), \|\varphi\|_{L^\infty(W)} \le 1} \left| \E(\varphi(X_t) \, | \, t< T_W) - \int_W \varphi d\nu \right| \le C \exp(-(\lambda_2-\lambda_1)t),
\end{equation}
where, as explained above, $-\lambda_1 > - \lambda_2$ are the two first eigenvalues of the operator $L=-\nabla V\cdot \nabla + \beta^{-1} \Delta$ considered on $W$ with zero Dirichlet boundary conditions on $\partial W$. 
\end{proposition}
The proof of this proposition is based on a spectral decomposition of the solution to the associated Feynman-Kac partial differential equation, and the simplest proof actually requires the law of $X_0$ to have a density with respect to $\mu$ which is in $L^2_\mu$.

This Proposition has two consequences. First, if the process $X_t$ enters a state and stays sufficiently long in the state, then its marginal in time is close to the QSD. Second, one could also think of considering the following interacting particle system, which samples the law of $X_t$ conditionally to the fact that $t< T_W$:
\begin{itemize}
\item Consider $N$ initial conditions $(X_0^n)_{1 \le n \le N}$ distributed independently according to a given law with support in $W$;
\item Let them evolve according to~\eqref{eq:sde}, driven by independent Brownian motions;
\item Each time a replica leaves the state, it is killed, another replica is duplicated and the new walker then evolves again independently of the others.
\end{itemize}
This is the so-called Flemming-Viot process \cite{burdzy-holyst-march-00,grigorescu-kang-04,ferrari-maric-07,lobus-08}. In the limit of infinitely many replicas, the replicas are distributed according to the law of $X_t$ conditionally to the fact that $t< T_W$, so that, in the long-time limit, the replicas are distributed according to the QSD (according to Proposition~\ref{prop:longtime_QSD}).

Finally, another crucial property of the QSD is the following:
\begin{proposition}\label{prop:QSD_out}
Let $X_0$ be distributed according to the QSD, and let us consider $(X_t)_{t \ge 0}$ solution to~\eqref{eq:sde}. Let us recall that $T_W=\inf\{t \ge 0, X_t \not \in W\}$. Then, 
\begin{itemize}
\item The law of the exit time $T_W$ is exponential with parameter $\lambda_1$;
\item The law of the exit point $X_{T_W}$ is $\displaystyle{\left(-\frac{1}{\beta \lambda_1} \frac{\partial u_1}{\partial n}\exp(-\beta V)\right) d\lambda_{\partial W}}$ where $n$ denotes the unit outward normal to $W$ and $\lambda_{\partial W}$ is the Lebesgue measure on $\partial W$;
\item The two random variables $T_W$ and $X_{T_W}$ are independent.
\end{itemize}
\end{proposition}
For proofs of these well-known properties of the QSD in our precise setting, we refer to~\cite{le-bris-lelievre-luskin-perez-11}.

In this context, one could state that the region $W$ is a metastable state for the dynamics~\eqref{eq:sde} if the typical time it takes to leave $W$ is large compared to the typical time it takes to reach the QSD (which is $1/(\lambda_2-\lambda_1)$ according to~\eqref{eq:cv_qsd}), namely:
\begin{equation}\label{eq:metastab_QSD}
\text{The probability } \P\left( T_W < \frac{1}{\lambda_2-\lambda_1}\right) \text{ is close to zero.}
\end{equation}
In some sense, this characterization~\eqref{eq:metastab_QSD} of metastability is the counterpart, in terms of QSD, of the two characterizations~\eqref{eq:metastab_LSI1} and~\eqref{eq:metastab_LSI2} we introduced above in terms of LSI.
The difficulty to completely formalize~\eqref{eq:metastab_QSD} is that the law of the initial condition $X_0 \in W$ should be defined to make precise the law of $T_W$. For example, if $X_0$ is distributed according to the QSD, $T_W$ is exponential with parameter $\lambda_1$ and thus $\P\left( T_W < \frac{1}{\lambda_2-\lambda_1}\right)=1-\exp(-\lambda_1/(\lambda_2-\lambda_1))$ so that $W$ is metastable if $\lambda_1 \ll \lambda_2 - \lambda_1$. As in the LSI case where~\eqref{eq:metastab_LSI2} could be used to define a good reaction coordinate, the characterization~\eqref{eq:metastab_QSD} of metastability in terms of the QSD could be useful to define what a good partition of the configurational space (namely a good function ${\mathcal S}$) is.

\subsection{A first application: analysis of the parallel replica dynamics}
\label{sec:ParRep}

The notion of QSD can be used to analyze an algorithm called the parallel replica dynamics which has been introduced by A.F. Voter in~\cite{voter-98}, and which is based on some Markovianity assumption, as will become clear below. The QSD is a way to quantify the error introduced by this Markovianity assumption, and to explain in which context the algorithm is efficient. The aim of the parallel replica dynamics is to generate very efficiently a process $(S_t)_{t \ge 0}$ with values in $\N$, and which is close to $({\mathcal S}(X_t))_{t \ge 0}$. Indeed, in many cases, one is not interested in the details of the dynamics of $X_t$: only the hopping events from one state to another are of interest. One important requirement is that the two stochastic processes  $(S_t)_{t \ge 0}$ and $({\mathcal S}(X_t))_{t \ge 0}$ should be close in terms of {\em the law of the trajectories} (not only the time marginals, for example, see also Remark~\ref{rem:olla} above).

Let us first describe the algorithm. This is a three stage algorithm. In the  {\em decorrelation step} a reference walker evolves according to~\eqref{eq:sde}, up to a time it stayed sufficiently long in the same state. More precisely, a time denoted $\tau_{\rm corr}$ is introduced, and  one proceeds to the next stage at a time $t_0$ if and only if ${\mathcal S}(X_t)=k_0$ is constant over the time interval $[t_0-\tau_{\rm corr},t_0]$. During all this step, $S_t$ is by definition ${\mathcal S}(X_t)$ so that no error is introduced. Let us assume that the decorrelation step has been successful, and let us proceed to the {\em dephasing step}. It consists in introducing $N$ replicas of the reference walker in the state $W_{k_0}$ and to let them evolve sufficiently long, conditionally to the fact that they do not leave the state, and to retain their final position. In other words, the dephasing step is basically a realization of the Fleming-Viot process introduced in the previous section. During this stage which is of course done in parallel, the process $S_t$ is not evolved. Finally, the speed up comes from the last step called the {\em parallel step}. It consists in letting all the walkers evolve independently from the initial conditions obtained in the previous dephasing step. Then the first escape event is detected, namely $$n_0=\argmin_{n \in \{1,\ldots,N\}} \{ T^n_W \}$$ where $T^n_W$ is the first escape time from $W$ for the $n$-th replica. Then, $S_t=k_0$ over the time interval $[t_0,t_0+ N T^{n_0}_W]$, and one proceeds to a new decorrelation step, the reference walker starting from the exit point $X^{n_0}_{T^{n_0}_W}$. The speed-up of course comes from the fact we consider only the first escape event among the $N$ walkers. As will become clear below, this event occurs in a time $N$ times smaller than the time it would take for a single walker to leave the state.

Let us now discuss the error analysis of this procedure. A first remark is that at the end of the decorrelation step, if $\tau_{\rm corr}$ has been chosen sufficiently large (namely $\displaystyle{\tau_{\rm corr} \gg \frac{1}{\lambda_2-\lambda_1}}$, according to~\eqref{eq:cv_qsd}), it is reasonable to assume that $X_{t_0}$ is approximately distributed according to the QSD. Thus, according to Proposition~\ref{prop:QSD_out}, the time it still remains to go out of the state $W_{k_0}$ is exponentially distributed, and independent of the exit point. Then, the aim of the dephasing step is clear: one wants to obtain $N$ initial conditions independently and identically distributed according to the QSD. The parallel step is thus fully justified. Indeed, concerning the time spent in the state $k_0$, since $(T^{1}_W,\ldots,T^N_W)$ are $N$ i.i.d. exponential random variables, $T^{n_0}_W=\min_n T^n_W$ is also an exponential random variable and $N T^{n_0}_W$ has the same law as $T^1_W$: this explains why the simulation clock is advanced by the amount of time $N T^{n_0}_W$ at the end of the parallel step. Moreover, concerning the exit point, since the exit time and the exit point are independent random variables when starting from the QSD, considering $X^{n_0}_{T_W^{n_0}}$ as the exit point is correct in terms of distribution (it has the same law as $X^1_{T_W^1}$).

In summary, the crucial parameter is $\tau_{\rm corr}$, which is used in the decorrelation step. The error which is made by one iteration of the algorithm can be formalized by considering:
$$e(t)=\sup_{f:\R_+ \times \partial W \to \R, \|f\|_{L^\infty} \le 1} \left|\E(f(T_W-t,X_{T_W})|T_W \ge t) - \E_\nu(f(T_W,X_{T_W})) \right|,$$
where $\E_\nu$ here denotes an expectation over functionals of $X_t$, the initial condition $X_0$ being distributed according to $\nu$. In words, $e(\tau_{\rm corr})$ measures the difference of what would have been the law of the couple of random variables (exit time, exit point) if the simulation of the reference walker would have been continued, compared to the law of the same couple of random variables, if we assume that the reference walker is distributed according to the QSD. A slight adaptation of Proposition~\ref{prop:longtime_QSD} above shows that $e(t) \le C \exp(-(\lambda_2 - \lambda_1)t)$ so that $\tau_{\rm corr}$ should be chosen larger than $1/(\lambda_2 - \lambda_1)$. The constant $C$ here depends on the law of the initial condition in the state. On the other hand, $\tau_{\rm corr}$ should be chosen smaller than the typical time it takes to leave the state, in order for the decorrelation step to have a chance to be successful. With these two requirements on $\tau_{\rm corr}$, it thus becomes clear that this algorithm is efficient if most of the states are metastable, in the sense of~\eqref{eq:metastab_QSD}, otherwise the decorrelation step will never be successful.

In conclusion, the interest of the QSD in this context is twofold: (i) it enables to understand how large $\tau_{\rm corr}$ should be in order not to introduce too much error in one iteration of the algorithm and (ii) it helps to define the assumptions required for the algorithm to be efficient.

\subsection{A second application: going from continuous state space dynamics to kinetic Monte Carlo models}
\label{sec:kMC}

The notion of QSD could also be useful in order to formalize the construction of discrete state space Markov models (so called kinetic Monte Carlo models~\cite{voter-05} in the context of molecular dynamics) from continuous state space Markov models such as~\eqref{eq:sde}. Let us recall that a stochastic process $S_t$ with values in $\N$ is Markovian if and only if: (M1) the list of visited states (forgetting about time) is Markovian and (M2) once $S_t$ takes a new value (it enters a new state), the time it takes to leave this state is exponentially distributed and independent of the next visited state.

There are basically two approaches to build such a connection. In the so-called milestoning approach~\cite{schuette-noe-lu-sarich-vanden-einjden-11,faradjian-elber-04}, one considers some disjoint subsets (the milestones) of the configuration space (think of small balls around the local minima of the potential $V$) and one considers the last milestone visited by $X_t$. The interest of this approach is that in the limit of very small subsets, the first requirement above (M1) is naturally satisfied. On the other hand, satisfying (M2) is more involved, and requires some assumptions related to the metastability of the process. Here metastability basically means that the time spent outside the milestones is very small compared to the time spent in the milestones, see~\cite{bovier-eckhoff-gayrard-klein-04,bovier-gayrard-klein-05}. The main drawback of this approach is that if the milestones are too small, the process $X_t$ spends a significant time outside of the milestones, so that the stochastic process built as ``the last visited milestone'' may not be a sufficiently fine coarse-grained description of the original process $X_t$, in order to extract useful macroscopic information (change of conformation of a molecular system, for example). Of course, this depends a lot on the system at hand.

The second natural approach which has been followed in the previous section and by many authors~\cite{huisinga-meyn-schuette-04,sarich-noe-schuette-10,voter-05} is to consider a full partition of the state space, and to consider at a given time $t$, in which state the process is. In the previous notation, it thus consists in considering~${\mathcal S}(X_t)$. Again, the process ${\mathcal S}(X_t)$ has no reason to be Markovian. The approximation suggested by the approach outlined above is to introduce a Markov process $S_t$ as follows: using the notation of the previous section, when $S_t$ jumps to a new value $k_0$ (one can imagine that the underlying process $X_t$ just enters a new state $W_{k_0}$), the time it takes to leave the value $k_0$ is exponentially distributed with parameter $\lambda_1$, and, independently, the next visited state is drawn according to the exit point distribution $\displaystyle{\left(-\frac{1}{\beta \lambda_1} \frac{\partial u_1}{\partial n}\exp(-\beta V)\right) d\lambda_{\partial W}}$. Here, we used the notation of the previous section: in particular, $(\lambda_1,u_1)$ are the first eigenvalue and eigenfunction of the operator $L$ on $W$, with zero Dirichlet boundary condition on $\partial W$. This procedure would be exact if, as soon as $X_t$ entered a new state, it would immediately be distributed according to the QSD. The error introduced by this coarse-grained description is thus related to the metastability of the original process, namely to the fact that when it enters a new state, it reaches the QSD before leaving the state. Contrary to the previous results presented in this paper, we have not yet fully formalized these ideas from a mathematical viewpoint. Notice that in the approach we propose here, the kernel of the approximating Markov process is computed using the QSD as an initial distribution in a given state, in contrast to what can be found usually in the literature, namely starting from the canonical measure $\mu$ restricted to the state. The interest of starting with the QSD is that the underlying assumptions ruling a Markov process (see (M1) and (M2) above) are automatically satisfied.

\begin{acknowledgement} I would like to thank my co-authors on these subjects (Chris Chipot, Nicolas Chopin, Giovanni Ciccotti, Brad Dickson, Benjamin Jourdain, Claude Le Bris, Fr\'ed\'eric Legoll, Mitch Luskin, Kimiya Minoukadeh, Stefano Olla, Danny Perrez, Mathias Rousset, Raphael Roux, Gabriel Stoltz and Eric Vanden-Einjden) as well as F\'elix Otto who introduced me to the so-called two-scale analysis for logarithmic Sobolev inequalities and Art Voter for very useful and inspiring discussions. This work is supported by the Agence Nationale de la Recherche, under grant ANR-09-BLAN-0216-01 (MEGAS).
\end{acknowledgement}


\end{document}